\documentclass[11pt]{article}
\usepackage[T1]{fontenc}
\usepackage[utf8]{inputenc}
\usepackage{authblk}

\usepackage[margin=1in]{geometry}  
\usepackage{graphicx}              
\usepackage{amsmath}               
\usepackage{amsfonts}              
\usepackage{amsthm}                

 \newtheorem{theorem}{Theorem}[section]

 \newtheorem{definition}[theorem]{Definition}
 
\theoremstyle{remark}

\title{A note on impossibility of uniformly non-oscillatory approximation \thanks{Carried out work is
    financially supported by DST India through project \#
    SR/FTP/MS-015/2011}
}
\author{Ritesh Kumar Dubey\footnote{\texttt{mail-to:
      riteshkd@gmail.com; riteshkumar.d@res.srmuniv.ac.in}}\\
Research Institute, SRM University, Tamilnadu\\ India}

\begin{document}
\maketitle

\begin{abstract}
In this note we show that it is impossible to have data independent
non-oscillatory three point finite difference scheme irrespective of its accuracy for a
scalar hyperbolic initial value problem.
\end{abstract}

{\bf keywords} Induced oscillations and numerical stability; affine
  combination; smoothness parameter; finite difference
  schemes; transport equation.

\section{Introduction}\label{sec1}
It is well elaborated in literature that evolution of discontinuities
are inevitable in the solution of scalar hyperbolic conservation laws
and artificial induced oscillations may appear in a numerical
approximation in the vicinity of discontinuity
\cite{Leveque,Laney,Toro}. One of the main emphasis in the numerical
approximation of hyperbolic conservation laws is to construct schemes
which are capable of yielding non-oscillatory approximation
\cite{Harten2,Shu1988,nessyahu} see also \cite{Shu1998}. The numerical
oscillations by high order schemes are well understood e.g., using
modified equation analysis. In \cite{Lax2006}, Lax showed that induced
oscillations analogous to the Gibbs phenomena must be present when the
solution is approximated by a difference scheme that is more than
first order accurate. Recently the cause of induced oscillations in
the solution even by first order monotone finite difference scheme for
hyperbolic problem is investigated in
\cite{Lefloach,Breuss2,Jiequan,Jiequan2}.
This work starts with one motivational example to show that
induced oscillations depend on the initial data. We define the
notion of data dependent stability and uniformly non-oscillatory
approximation. We finally show that it is impossible to have uniformly
non-oscillatory approximation by any three point scheme irrespective
of its accuracy. 
\section{Uniformly non-oscillatory approximation}
We consider the following simple linear initial value problem, 
\begin{equation}
\frac{\partial u(x,t)}{\partial t} + a\frac{\partial u(x,t)}{\partial x}=0, a\neq0, u(x,0)=u_{0}(x) \label{transporteq}
\end{equation}
where $u(x,t)$ is a scalar function of space $x$ and time $t$
variables, and $u_{0}(x)$ is a piece-wise smooth function with possible
discontinuities, and characteristic function $a(x,t)$ is a smooth
function. We discretize the space and time
variable in to a computational mesh using spatial and temporal step
$h$ and $k$ respectively with grid points $(x_{i}=ih,t_n=nk)$. Let
$\lambda=\frac{k}{h}$ and $u_{i}^{n}$ represents and approximation of
the value $u(ih,nk)$.
\subsection{Lax-Wendroff (LxW) scheme}
\begin{figure}[!htb]
  \begin{tabular}{cc}
    \includegraphics[scale=0.35]{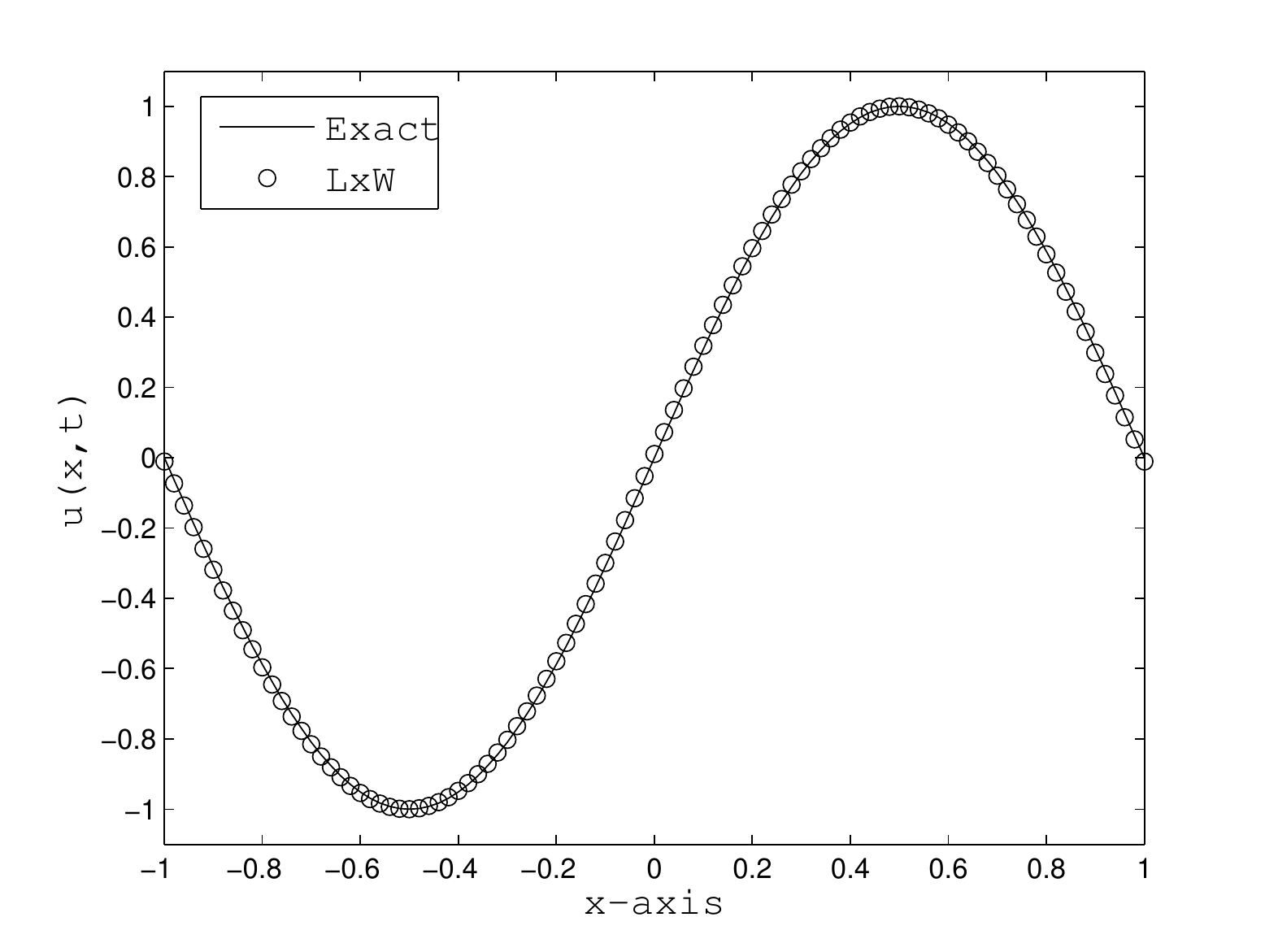}&  \includegraphics[scale=0.35]{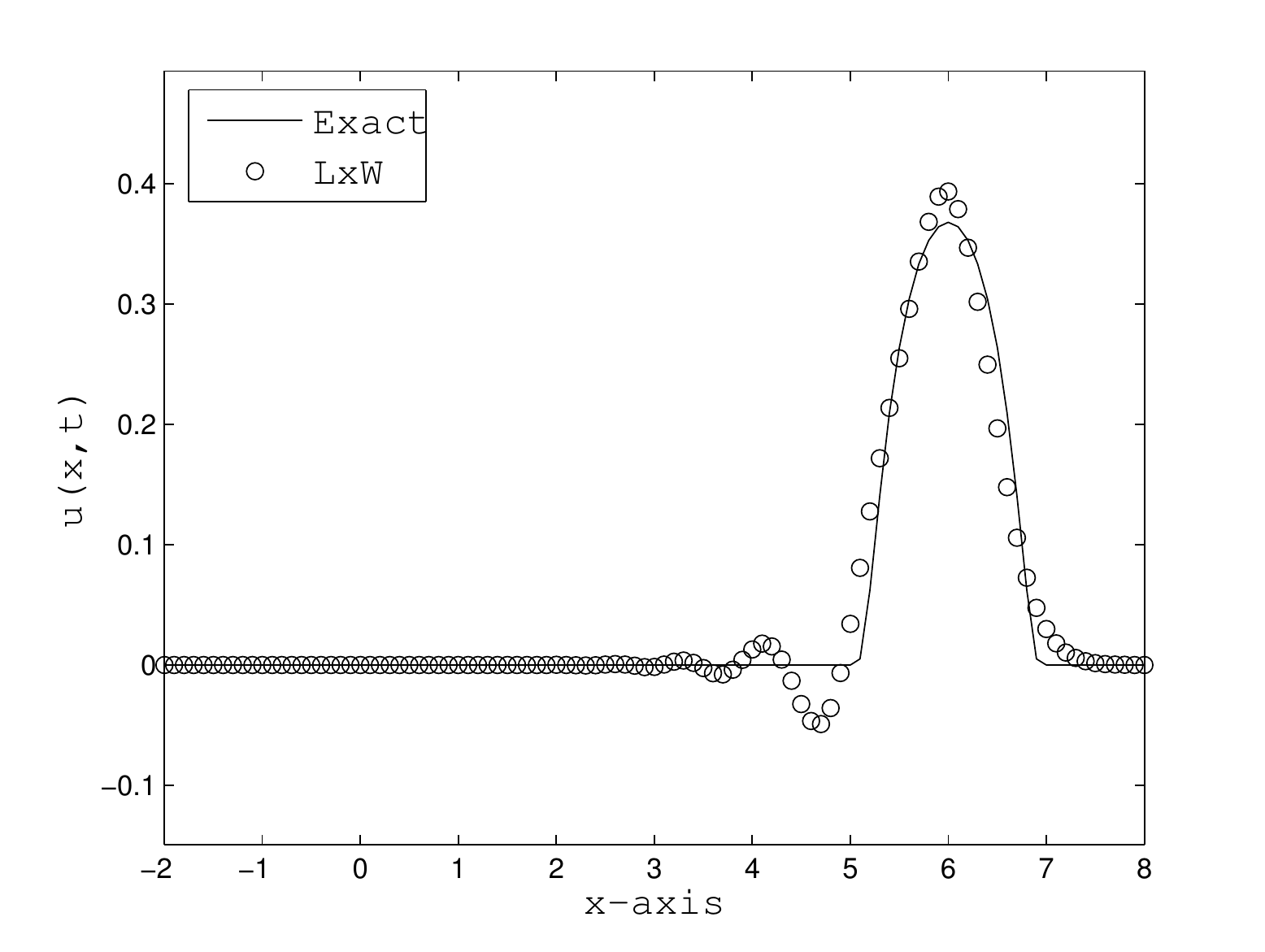}\\
    (a) & (b)
  \end{tabular}
  \caption{Solution using data $CFL=0.8, t=6, N=100$. Both initial condition are with out discontinuity (a) No induce oscillation (b) induced oscillations\label{Fig1}}
\end{figure}

Consider LxW scheme for (\ref{transporteq}). It is second order
accurate therefore introduces oscillations in computed solution for
discontinuous initial condition and is well understood by
\cite{Lax2006}. We now consider, (\ref{transporteq}) corresponding to
smooth initial condition
\begin{itemize}
  \item[a] $u_{0}(x)=\sin(\pi\,x),\; x\in [-1:1] $
  \item[b] $u_{0}(x)=\left\{\begin{array}{cc} \exp^{\frac{-1}{1-x^2}} & if x\in [-1:1]\\ 0 & else \end{array}\right. $
\end{itemize}
The numerical results by LxW are given in Figure \ref{Fig1}. It is
clear from the results that even though both initial conditions are
smooth, LxW does not introduces oscillations in solution corresponding
to initial condition (a) while oscillations appears in the solution
for initial data (b). Also as mentioned above even a first order
monotone Lax-Friedrichs scheme exhibits local oscillations for
discontinuous solution see \cite{Breuss2}. Therefore it can be
concluded that the induced oscillations depends on the initial data
and can induce in to a numerical solution irrespective of its
accuracy. This makes it reasonable to analyze the oscillatory behavior
of any scheme with respect to the initial data.
\subsection{A non-oscillatory condition for two point}
The idea is to use the following simple concept from
geometry: Let each $P_{1}, P_{2}$ denotes two points in a finite
affine plane then there affine combination is defined by
\begin{equation}
P = P_{1} + \delta(P_{2}-P_{1})\label{affinetwopts1}
\end{equation}
Geometrically, $P$ represents a point on the line $L$ which passes
through $P_{1}$ and $P_{2}$, and if $0\leq \delta\leq1$, $P$ lies on
the line segment joining point $P_{1}$ and $P_{2}$. Then following the
method of characteristics and affine combination
(\ref{affinetwopts1}), an initial data independent non-oscillatory
solution of (\ref{transporteq}) can be ensured by an explicit two
point consistent difference schemes scheme of the form
\begin{equation}
u^{n+1}_{j} =\left\{\begin{array}{cc} u^{n}_{j} -D \Delta u^{n}_{j}, & \; \mbox{if} \; a(x,t)> 0,\\ 
u^{n}_{j} + D  \Delta_{+}u^{n}_{j}, & \; \mbox{if} \; a(x,t)< 0,\end{array}\right.\label{twopts1}
\end{equation}
provided $0 \leq D \leq1$. The coefficient $D$ depends on $a,\;
\lambda$. A simple Von-Neumann stability analysis also shows that
(\ref{twopts1}) is stable provided $D \leq1$.
We call scheme (\ref{twopts1}) uniformly non-oscillatory (UNO).
\section{Non-oscillatory condition for three point scheme}
In this section we show that it is impossible to have a initial data
independent uniformly non-oscillatory scheme with higher then two
points.  In \cite{Lax2006}, Lax pointed out the use of the following
smoothness parameter to devise hybrid scheme which can avoid the Gibbs
phenomena. It is defined similar to the one in
\cite{vanLeer1974,Sweby,Rkd1}
\begin{equation}
\theta_{i}^{n}=  \left\{\begin{array}{ccc}\displaystyle
\theta_{i}^{n,+}=&\frac{\Delta_{-}u_{i}^{n}}{\Delta_{+}u_{i}^{n}}&\;\; if\;\; a\geq
0\\ &\\\theta_{i}^{n,-}= &  \displaystyle\frac{\Delta_{+}u_{i}^{n}}{\Delta_{-}u_{i}^{n}}& \;\;if\;\;
 a <0.\end{array}\right.\label{deftheta}
\end{equation}
Using above affine combination we define
\begin{definition}\label{def1}
Consider a consistent three point scheme written in the form
\begin{equation}
u^{n+1}_{j} =\left\{\begin{array}{cc}
u^{n}_{j} + D(a\lambda;\theta_{i}) \Delta_{-} u^{n}_{j} & \; \mbox{if} \; a(x,t)> 0\\ 
u^{n}_{j} + D(a\lambda;\theta_{i}) \Delta_{+} u^{n}_{j} & \; \mbox{if} \; a(x,t)< 0,\end{array}\right. \label{defDDS}
\end{equation}
where $i\in\{j, j\pm1 \}$ and the parameter coefficient $D$ depends on
CFL number $a\lambda$ as well smoothness parameter
(\ref{deftheta}). The scheme (\ref{defDDS}) is uniformly non-oscillatory (UNO) if
\begin{equation}
  0 \leq  D(\theta_{i})\leq 1,\;  \forall\;\; \theta_{i}\in\mathbb{R}, \label{ddscondition}
\end{equation}
In case (\ref{ddscondition}) holds only for $\theta_{i}\in
\mathbb{R}\setminus \mathbb{S}$ where $\mathbb{S}\subset \mathbb{R}$,
scheme (\ref{defDDS}) is called Data dependent stable (DDS). Note that
three point scheme (\ref{defDDS}) using a condition on $\theta$
ensures that the updated value $u_{i}^{n+1}$ satisfies the following
local maximum principle
$$ \begin{array}{cc} u_{i-1}^{n}  \leq u_{i}^{n+1} \leq u_{i}^{n},&\; \mbox{if}\; a>0,\; \\
   u_{i}^{n}  \leq u_{i}^{n+1} \leq u_{i+1}^{n},&\; \mbox{if}\; a<0.\end{array}$$
\end{definition}

\subsection{Three point upwind linear schemes}
A generic numerical flux function for any linear three point upwind
scheme for (\ref{transporteq}) can be given by,
\begin{equation}
F_{i+\frac{1}{2}}= \left\{\begin{array}{cc}\alpha u^{n}_{i} +\beta u^{n}_{i-1}
&\;\; if\;\; a> 0\\\alpha u^{n}_{i+1} +\beta u^{n}_{i+2} &\;\; if\;\; a<0.\end{array}\right. \label{fluxupwind}
\end{equation} 
where for consistency the coefficients must satisfy $\alpha+\beta =a.$\footnote{A consistent
discretization requires i.e., $F(u,u)=au$}.
\begin{itemize}
  \item For $a> 0$, the resulting three point conservative
difference scheme  using flux (\ref{fluxupwind}) is
\begin{equation}
u_{i}^{n+1}=u_{i}^{n}-\lambda\left(\alpha \Delta_{-}u_{i}^{n} - \beta \Delta_{-}u_{i-1}^{n}\right) \label{tpus1}
\end{equation}
or
\begin{equation}
u_{i}^{n+1}=u_{i}^{n}-\lambda\left(\alpha  - \beta \frac{\Delta_{-}u_{i-1}^{n}}{\Delta_{-}u_{i}^{n}}\right) \Delta_{-}u_{i}^{n}\label{upwind1}
\end{equation}
By Definition \ref{def1}, approximation (\ref{upwind1}) is DDS provided 
$$0\leq \lambda(\alpha-\beta \theta_{i-1}^n)\leq 1.$$
or
$$-\alpha \leq \beta \theta_{i-1}^{n,+}\leq \frac{1-\alpha \lambda}{\lambda}.$$
which satisfies
\begin{subequations}
\begin{equation}
  -\frac{\alpha}{\beta} \leq \theta^{n,+}_{i-1}\leq \frac{1-\lambda\alpha}{\lambda \beta},\; \mbox{if}\; \beta\geq 0 \label{tpUSInq1a}
\end{equation}
or
\begin{equation}
  \frac{1-\lambda\alpha}{\lambda \beta}\leq \theta^{n,+}_{i-1}\leq   -\frac{\alpha}{\beta},\; \mbox{if}\; \beta\leq 0 \label{tpUSInq1b}
\end{equation}\label{tpUSInq1}
\end{subequations}
\item For $a<0$, resulting scheme using flux in (\ref{fluxupwind}) can be written as
  \begin{equation}
    u_{j}^{n+1} =u_{j}^{n}+ \lambda\left[\alpha -\beta \frac{\Delta_{+} u_{j+1}^{n}}{\Delta_{+} u_{j}^{n}} \right] \Delta_{+} u_{j}^{n}    
  \end{equation}
  which by similar calculation is DDS provided

  \begin{subequations}
    \label{tpUSInq2}
    \begin{equation}
      -\frac{\alpha}{\beta} \leq \theta^{n,-}_{j+1} \leq \frac{1- \lambda\alpha}{\lambda\beta},\; \mbox{if}\; \beta \geq 0  \label{tpUSInq2a}
    \end{equation}
    \begin{equation}
      \frac{1- \lambda\alpha}{\lambda\beta} \leq \theta^{n,-}_{j+1} \leq -\frac{\alpha}{\beta},\; \mbox{if}\; \beta \leq 0  \label{tpUSInq2b}
    \end{equation}
  \end{subequations}
\end{itemize}  
Note that DDS conditions (\ref{tpUSInq1}) and (\ref{tpUSInq2}) gives
conditions on the initial data in terms of smoothness parameter
$\theta$ such that local oscillations does not introduced by a three
point upwind schemes.  In Table \ref{tab:1}, Data dependent stability region in terms of $\theta$ is
given for classical upwind schemes.
\begin{table}
  \hspace{-1.5cm}\begin{tabular}{|c|c|c|c|c|c|}
    \hline
    scheme & Order & $\alpha$ & $\beta$ & CFL condition & DDS bound \\
    \hline
    Two point upwind & First & $a$ & $0$  & $0<a\lambda \leq 1$ & UNO \\
    \hline
    Three point upwind &Second&  $\frac{3a}{2}>0$ & $\frac{a}{2}>0$  &$0<a\lambda \leq \frac{1}{2}$ & $\theta^{n,+}_{i-1} \in  \left[-\frac{2-3a\lambda}{a\lambda},3\right]$ \\
    \hline
    \hline
    Beam-Warming &Second& $\frac{a}{2}(3-\lambda a)>0$ & $\frac{a}{2}(1-\lambda a)>0$ &$0<a\lambda <1$ & $\theta_{i-1}^{n,+} \in \left[-\frac{2-a\lambda(3-a\lambda)}{a\lambda(1-a\lambda)},\frac{3-a\lambda}{1-a\lambda}\right]$\\
    \hline Beam-Warming &Second& $\frac{a}{2}(3-\lambda a)>0$ & $\frac{a}{2}(1-\lambda a)<0$ &$1<a\lambda \leq 2$ & $\theta_{i-1}^{n,+} \in \left[-\frac{2-a\lambda(3-a\lambda)}{a\lambda(a\lambda-1)},\frac{3-a\lambda}{a\lambda-1}\right]$ \\
    \hline
  \end{tabular}      
\caption{Data dependent stability bounds for three point upwind schemes for $a>0$. Beam-Warming scheme changes its DDS interval with respect to $CFL$.}
\label{tab:1}       
\end{table}
\subsection{Three point centred linear schemes}
Consider the following generic consistent numerical flux function of three
point centred linear scheme for (\ref{transporteq})
\begin{equation}
F_{i+\frac{1}{2}}= \alpha u^{n}_{i+1} + \beta u^{n}_{i},\label{fluxcenterd}
\end{equation}
where again for consistency $\alpha+\beta =a$.
\begin{itemize}
  \item $a>0$ The resulting conservative approximation can be written
    as
\begin{equation}
u^{n+1}_{i} = u^{n}_{i}-\lambda [\alpha \Delta_{+}u^{n}_{i} + \beta
  \Delta{-}u^{n}_{i}] \label{tpcs1}
\end{equation}
which using (\ref{deftheta}) reduces to
\begin{equation}
  u^{n+1}_{i} =u^{n}_{i} -\lambda \left( \alpha \theta_{i}^{n,-}+\beta \right) \Delta_{-}u^{n}_{i},\; a>0\label{centred1}
\end{equation}
by Definition \ref{def1}, (\ref{centred1}) is non-oscillatory stable provided
\begin{equation}
-\beta\leq \alpha \theta_{i}^{n,-}\leq \frac{1-\beta \lambda}{\lambda}
\end{equation}
it reduces to
\begin{subequations}
  \begin{equation}
    -\frac{\beta}{\alpha} \leq \theta_{i}^{n,-} \leq \frac{1-\lambda \beta}{\lambda \alpha},\; \alpha>0
  \end{equation}
  \begin{equation}
    \frac{1-\lambda \beta}{\lambda \alpha} \leq \theta_{i}^{n,-} \leq -\frac{\beta}{\alpha} ,\; \alpha>0
  \end{equation}
\end{subequations}
which on inversion gives condition for (\ref{tpcs1}) to be DDS
\begin{subequations}
  \label{tpCSInq1}
\begin{equation}
\theta_{i}^{n,+} \in \left(-\infty, -\frac{\alpha}{\beta}\right]\cup
  \left[\frac{\alpha \lambda}{1-\beta \lambda},\infty
  \right),\; \alpha>0 \label{tpCSInqa}
\end{equation}
\begin{equation}
\theta_{i}^{n,+} \in \left(-\infty, \frac{\lambda \alpha}{1-\lambda \beta}\right]\cup \left[ -\frac{\alpha}{\beta},\infty\right), \; \alpha<0.\label{tpCSInqb}
\end{equation}
\end{subequations}
\item $a<0$ scheme (\ref{tpcs1}) can be written as,
\begin{equation}
  u^{n+1}_{i} = u^{n}_{i} +\lambda \left(\alpha +\beta \theta_{i}^{n,+} \right) \Delta_{+}u^{n}_{i},\; a<0 \label{centred2}  
\end{equation}
which using similar calculation is DDS provided
\begin{subequations}
  \label{tpCSInq2}
  \begin{equation}
    \theta_{i}^{n,-} \in \left(-\infty, -\frac{\beta}{\alpha}\right]\cup\left[\frac{\lambda \beta}{1-\lambda \alpha}, \infty\right),\;\mbox{if}\; \beta>0 \label{tpCSInq2a}
  \end{equation}
  \begin{equation}
    \theta_{i}^{n,-} \in \left(-\infty, \frac{\lambda \beta}{1-\lambda \alpha}\right]\cup\left[-\frac{\beta}{\alpha}, \infty\right),\;\mbox{if}\; \beta<0 \label{tpCSInq2b}
  \end{equation}
\end{subequations}
\end{itemize}
The DDS conditions (\ref{tpCSInq1}) and (\ref{tpCSInq2}) classify the
data type in terms of the smoothness parameter $\theta$ such that
local oscillations does not get introduced by three point centred
schemes. In Table \ref{tab:2}, data dependent stability reason is
given for classical three point centred schemes which justifies the
numerical oscillations by these schemes irrespective of their accuracy.
\begin{table}
  \hspace{-1.5cm}\begin{tabular}{|c|c|c|c|c|c|}
    \hline scheme & Order &CFL & $\alpha$ & $\beta$ & DDS bound
    \\ \hline Lax-Friedrichs & First &$0<a\lambda \leq1$ &
    $\frac{1}{2\lambda}(a\lambda -1)<0$ & $\frac{1}{2\lambda}(a\lambda
    +1)>0$ & $ \theta^{n,+}_{i} \in (\infty,
    -1]\cup[\frac{1-a\lambda}{1+a\lambda},\infty)] $\\ \hline FTCS &
  Second & $0<a\lambda \leq1$ &$ \frac{a}{2}>0$ & $ \frac{a}{2}>0$ &
  $ \theta^{n,+}_{i} \in
  [-\infty,-1]\cup\left[\frac{a\lambda}{2-a\lambda},\infty\right)$\\ \hline
    Lax-Wendroff& Second&$0<a\lambda \leq1$ & $
    \frac{a}{2}(1-a\lambda)>0$ &$ \frac{a}{2}(1+a \lambda)>0$ &
    $\theta^{n,+}_{i} \in \left(-\infty,
    -\frac{1-a\lambda}{1+a\lambda}\right]\cup\left[\frac{a\lambda}{2+a\lambda}\right),\infty)$\\ \hline
  \end{tabular}      
  \caption{Data dependent stability bounds for three point centred schemes under CFL $|a|\lambda \leq 1$}
  \label{tab:2}       
\end{table}
In order to achieve uniformly non-oscillatory (UNO) approximation by
an upwind or centred scheme one needs to choose coefficients
$\alpha,\, \beta$ in numerical flux (\ref{fluxupwind}) or
(\ref{fluxcenterd}) such that resulting scheme is DDS $\forall\;
\theta_{i\mp1}^{n,\pm}\, \in \mathbb{R}$. Bounds in (\ref{tpUSInq1})
and (\ref{tpUSInq2}) or (\ref{tpCSInq1}) and (\ref{tpCSInq2}) shows
that scheme three point upwind or centred are not UNO and will always
introduce local numerical oscillation except for the choice
$\alpha=a,\; \beta=0$ in (\ref{fluxupwind}) or the choice $\alpha=0,\;
a>0$ and $\beta=0,\; a<0$ in (\ref{fluxcenterd}). This along with
consistency requirement results in to {\it classical first order upwind
  scheme}.\\
\subsubsection{ Example: A hybrid second order UNO scheme}
The second order essentially non-oscillatory {ENO} scheme \cite{Shu1998} for transport problem is UNO
under the CFL condition $0\leq \lambda a \leq \frac{1}{2}$.\\
Proof: Let $a>0$, then the ENO second order reconstruction is
$$ \hat{F}_{i+\frac{1}{2}}= \left\{\begin{array}{cc}\displaystyle
-\frac{a}{2}u_{i-1}+\frac{3a}{2}u_i &\; \mbox{if}\;
|u_{i-1}-u_i|<|u_{i}-u_{i+1}|\\ &\\ \displaystyle
\frac{a}{2}u_{i}+\frac{a}{2}u_{i+1} &
\mbox{else}\end{array}\right.$$
The resulting ENO scheme for the transport problem (\ref{transporteq}) can be written as
\begin{equation}
{u}_{i}^{n+1}=u_{i}-D\Delta_{-}{u}_{i}   
\end{equation}
where 
$$
D=\left\{\begin{array}{ll}
a\lambda \left[-\frac{1}{2}\frac{\Delta{-}u_{i-1}}{\Delta_{-}u_{i}}+\frac{3}{2}\right] & \mbox{if}\; |\Delta_-{u_{i}}|<|\Delta_+{u_{i}}|\;\mbox{and}\;|\Delta_-{u_{i-1}}|<|\Delta_-{u_{i}}|
\\a\lambda  & \mbox{if}\; |\Delta_-{u_{i}}|<|\Delta_+{u_{i}}|\; \mbox{and}\;|\Delta_-{u_{i-1}}|>|\Delta_-{u_{i}}| 
\\a\lambda  \left[\frac{1}{2}(\frac{\Delta_+u_{i}}{\Delta_-u_{i}}-\frac{\Delta_-u_{i-1}}{\Delta_-u_{i}})+1\right] & \mbox{if}\; |\Delta_-{u_{i}}|>|\Delta_+{u_{i}}|\;\mbox{and}\; |\Delta_-{u_{i-1}}|<|\Delta_-{u_{i}}|
\\ a\lambda \left[\frac{1}{2}\frac{\Delta_+u_{i}}{\Delta_-u_{i}}+\frac{1}{2}\right] & \mbox{if}\;|\Delta_-{u_{i}}|>|\Delta_+{u_{i}}|\;\mbox{and}\;|\Delta_-{u_{i-1}}|>|\Delta_-{u_{i}}|
\end{array}\right.
$$ Note under the CFL condition $0< \lambda a \leq\frac{1}{2}$, $0
\leq D \leq 1$, thus scheme is UNO.
\section{Conclusion}
An approach to find non-oscillatory condition on initial data in terms
of smoothness parameter is demonstrated for three point schemes. It
shows that it is impossible to have a data independent non-oscillatory
three point scheme. In future it will be inter sting to extend this
approach to analyze ENO schemes for their non-oscillatory condition
and construction of hybrid high order UNO scheme.


\begin{thebibliography}{}
%
%

  \bibitem{Lax2006} Lax P. D., Gibbs Phenomena, Journal of Scientific
    Computing, Volume 28, Issue 2-3, pp 445-449, 2006.
\bibitem{Leveque} R. J. LeVeque, Numerical methods for conservation
  laws, Lectures in Mathematics ETH Zurich, 2nd ed. Basel,
  Birkh\"{a}user Verlag, 1992.
\bibitem{Laney} C B Laney, Computational Gas-dynamics, Cambridge
  University press, 1998.
\bibitem{Lefloach} P G Lefloach and J G Liu, Generalized monotone
  schemes, discrete paths of extrema and discrete entropy conditions,
  Math. Comp. 68,pp 1025-1055, 1999.
\bibitem{Breuss2} M. Breu\ss, An Analysis of the Influence of Data
  Extrema on some first and second order central approximations of
  Hyperbolic Conservation Laws, ESAIM: Mathematical Modelling and
  Numerical Analysis, Volume 39, No. 5, 965-994,2005.
\bibitem{Jiequan} Jiequan Li, H. Tang, G. Warneke and L. Zhang, Local
  Oscillations in finite difference solutions of hyperbolic
  conservation laws, Mathematics of Computation, 78, 1997-2018,2009.
\bibitem{Jiequan2} Jiequan Li, Z. Yang, Heuristic modified equation
  analysis of oscillations in numerical solutions of conservation
  laws, SIAM J. Num. Analysis, 49, 2386-2406,2011.
\bibitem{Harten2} Harten A. and Osher S., "Uniformly high-order
  accurate non-oscillatory schemes I," SIAM J. Numer. Anal.,
  v. 24, pp. 279-309, 1987.
\bibitem{Shu1988}  Chi-Wang Shu, Stanley Osher,
  Efficient implementation of essentially non-oscillatory shock-capturing schemes,
  Journal of Computational Physics, Volume 77, Issue 2, Pages 439-471,1988.
\bibitem{Shu1998} Chi-Wang Shu, Essentially non-oscillatory and
  weighted essentially non-oscillatory schemes for hyperbolic
  conservation laws, Advanced Numerical Approximation of Nonlinear
  Hyperbolic Equations, Lecture Notes in Mathematics Volume 1697, pp
  325-432, 1998.  
\bibitem{nessyahu} H. Nessyahu and E. Tadmor, Non-Oscillatory central
  differencing for hyperbolic conservation laws,
  J. Comp. Phys. 87, 408-436, 1990.
\bibitem{Sweby} Sweby P. K., High resolution schemes using flux
  limiters for hyperbolic conservation laws, SIAM
  J. Numer. Anal. 21, 995-1011, 1984.
\bibitem{Toro} Toro E. F., Riemann solvers and numerical methods for
  fluid dynamics, A practical introduction 2nd edition, Springer,
  1999.
\bibitem{Rkd1} Dubey Ritesh Kumar, Flux limited schemes: Their
  classification and accuracy based on total variation stability
  regions, Appl. Math. Comput. 224, pp 325-336, 2013. 
\bibitem{Yongi}
  Yongqi Wang, Kolumban Hutter, \newblock{Comparisons of numerical methods with respect to
    convectively dominated problems},
  \newblock{Int. J. Numer. Meth. Fluids} 37,721-745, 2001.
\bibitem{Ewing} R.E. Ewing, Hong Wang,
  \newblock{A summary of numerical methods for time-dependent advection-
    dominated partial differential equations},
  \newblock{J. Comput. Appl. Math.} 128, 423-445, 2001.
\bibitem{vanLeer1974}
Bram van Leer.
\newblock Towards the ultimate conservative difference scheme. ii. monotonicity
  and conservation combined in a second-order scheme.
\newblock {\em Journal of Computational Physics}, 14(4):361 -- 370, 1974.

\end{thebibliography}
\end{document}